\documentclass[11pt]{article}
\usepackage{bbm}
\usepackage{mathrsfs}
\usepackage{amsmath}
\usepackage{amsthm}
\usepackage{amsfonts}
\usepackage{amssymb}
\usepackage{latexsym}
\usepackage{amstext}
\usepackage{array}
\usepackage{graphicx}
\usepackage{MnSymbol}
\usepackage{bbm}
\usepackage{amsfonts}
\usepackage{amssymb}
\usepackage{ams,amsmath,amsthm}
\usepackage{multirow}
\usepackage{booktabs}

\newtheorem{thm}[subsection]{Theorem}
\newtheorem{cor}[subsection]{Corollary }
\newtheorem{Def}[subsection]{Definition}
\newtheorem{lem}[subsection]{Lemma}
\newtheorem{remark}[subsection]{Remark}

\newtheorem{prop}[subsection]{Proposition}
\newtheorem{exm}[subsection]{Example}
\newcommand{\bthm}{\begin{thm} }
\newcommand{\ethm}{\end{thm} }
\newcommand{\bpro}{\begin{prop}}
\newcommand{\epro}{\end{prop}}
\newcommand{\bdf}{\begin{Def}}
\newcommand{\edf}{\end{Def}}
\newcommand{\bexm}{\begin{exm}}
\newcommand{\eexm}{\end{exm}}
\newcommand{\blem}{\begin{lem}}
\newcommand{\elem}{\end{lem}}
\newcommand{\bpf}{\begin{proof}}
\newcommand{\epf}{\end{proof}}
\newcommand{\bcor}{\begin{cor}}
\newcommand{\ecor}{\end{cor}}
\newcommand{\beq}{\begin{equation}}
\newcommand{\eeq}{\end{equation}}
\newcommand{\ba}{\begin{array}}
\newcommand{\ea}{\end{array}}
\newcommand{\bea}{\begin{eqnarray}}
\newcommand{\eea}{\end{eqnarray}}

\newcommand{\brem}{\begin{remark}}
\newcommand{\erem}{\end{remark}}

\begin{document}


\title{Explicit primality criteria for $h\cdot 2^n\pm 1$}

\author{
Yingpu Deng and Dandan Huang\\\\
Key Laboratory of Mathematics Mechanization,\\
NCMIS, Academy of Mathematics and Systems Science,\\
Chinese Academy of Sciences, Beijing 100190, P.R. China\\
Email: \{dengyp, hdd\}@amss.ac.cn}

\date{}
\maketitle

\begin{abstract}
We describe an explicit generalized Lucasian test to determine the primality of
numbers $h\cdot2^n\pm1$ when $h\nequiv0\pmod{17}$. This test is by means of fixed seeds which depend only
on $h$. In particular when $h=16^m-1$ with $m$ odd,
our paper gives a primality test with some fixed seeds depending only on $h$. Comparing the results of W. Bosma(1993) and P.
Berrizbeitia and T. G. Berry(2004), our result adds new values of $h$ along with this line. Octic and bioctic reciprocity are used to deduce our result.

\end{abstract}

\section{Introduction}

In this paper we consider primality tests for integers $M$ of the form $h\cdot2^n\pm1$ with $h$ odd.
Primality tests for numbers of such form have been noticed since Lucas \cite{lu} and Lehmer \cite{le} gave the celebrated Lucas-Lehmer primality test for Mersenne numbers, using properties of the Lucas sequences. Here, we recall this famous primality test:

    \textbf{Lucas-Lehmer test.}\quad Let $M_p=2^p-1$ be Mersenne number, where $p$ is an odd prime. Define a sequence $\{u_k\}$ as follows: $u_0=4$ and $u_k=u_{k-1}^2-2$ for $k\geq1$. Then $M_p$ is a prime if and only if $u_{p-2}\equiv0\pmod{M_p}$.

We call a sequence $\{u_k\mid k\geq0\}$ is a Lucasian sequence if the recurrence relation is $u_k=u_{k-1}^2-2$ for $k\geq1$ and $u_0$ is called the seed of the sequence.
In 1993, Bosma \cite{bo} posed the problem whether there exists finitely many seeds depending only on $h$ of some Lucasian sequences for which the sequences can determine the primality of $h\cdot2^n\pm1$. In the same paper Bosma exhibited that a finite set of pairs $(d_k, \alpha_k)$ with $d_k\in\mathbb{Z}$ and
$\alpha_k\in\mathbb{Q}(\sqrt{d_k})$ always exists, such that, for any $n$, one of
the pairs determines the primality of $M=h\cdot2^n\pm1$ for $h<10^5$, except $h=4^m-1$ with $m>0$.

For $h\nequiv0\pmod{3}$, we know that, from the result in \cite{ri}(or see \cite{bo}), such a single seed of the Lucasian sequence exists. For $h\nequiv0\pmod{5}$, using biquadratic reciprocity, in \cite{be}, P.
Berrizbeitia and T. G. Berry have shown that there exists
a Lucasian sequence with a single seed independent of $n$ to test the primality of numbers of the form
$M=h\cdot2^n\pm1$. In particular, for $h=4^m-1$ with $m$ odd, this holds.

In this paper, we will prove that, for fixed $h\nequiv0\pmod{17}$, there are
some generalized Lucasian sequences with fixed seeds independent of $n$ which can
determine the primality  for integers of the form $M=h\cdot2^n\pm1$. In particular
our paper further shows that the fixed seeds exist for $h=16^m-1$ with $m$ odd. Octic and bioctic reciprocity are used to deduce our result, and the key point is that we use two or four sequences rather than a single sequence involved in the above mentioned works.

\smallskip

\section{Octic and Bioctic reciprocity}

What we state in this section can be found in \cite[Chapter 14]{ir} and \cite[Chapter 14]{ev}.

Let $\zeta_m=e^{2\pi \sqrt{-1}/m}$ be the complex primitive $m$-th root of unity, and let
$D=\mathbb{Z}[\zeta_m]$ be the ring of integers of the cyclotomic field $\mathbb{Q}(\zeta_m)$. Let $\mathfrak{p}$ be a prime ideal of $D$ lying over a rational prime $p$ with gcd$(p, m)=1$. For every $\alpha\in D$, the $m$-th power residue
symbol $\left(\frac{\alpha}{\mathfrak{p}}\right)_m$ is defined by:

$\mathrm{(1)}$ If $\alpha\in\mathfrak{p}$, then $\left(\frac{\alpha}{\mathfrak{p}}\right)_m$ = 0.

$\mathrm{(2)}$ If $\alpha\notin\mathfrak{p}$, then $\left(\frac{\alpha}{\mathfrak{p}}\right)_m$ = $\zeta_m^i$ with $i\in\mathbb{Z}$, where $\zeta_m^i$
is the unique $m$-th root of unity in $D$ such that
$$\alpha^{(\mbox{N}(\mathfrak{p})-1)/m}\equiv \zeta_m^i\quad \pmod{\mathfrak{p}},$$
where N$(\mathfrak{p})$ is the absolute norm of the ideal $\mathfrak{p}$.

$\mathrm{(3)}$ If $\mathfrak{a}\subset D$ is an arbitrary ideal and $\mathfrak{a}=\prod\mathfrak{p}_i^{n_i}$ is its factorization
 as a product of prime ideals, then
 $$\left(\frac{\alpha}{\mathfrak{a}}\right)_m=\prod\left(\frac{\alpha}{\mathfrak{p}_i}\right)_m^{n_i}.$$
We set $\left(\frac{\alpha}{D}\right)_m$ = 1.

 $\mathrm{(4)}$ If $\beta\in D$ and $\beta$ is prime to $m$ define
 $\left(\frac{\alpha}{\beta}\right)_m=\left(\frac{\alpha}{\beta D}\right)_m$.

Let $m=l^n$ ($\neq2, 4$), where $l$ is a prime and $n$ is a positive integer.
An element $\alpha\in D$ is said to be primary if $\alpha$ is coprime with $m$
and $\epsilon_c(\alpha)=(-1)^M$, where $M=(\mbox{N}(\alpha D)-1)/m$, $c$ is some given
integer, and $\epsilon_c(\alpha)$ is a power of $\zeta_m$ whose explicit definition
can be found in \cite[Chapter 14]{ev}. The following theorem can be found in \cite[Chapter 14, Th.14.3.1, p. 474]{ev}.

\smallskip

\begin{thm}(Eisenstein's Reciprocity Law)\label{Eis}
Let $m=l^n$ ($\neq2, 4$), where $l$ is a prime and $n$ is a positive integer.
 Let $a$ be a rational prime with $gcd(a, m)=1$, and let $\alpha$ be a primary integer of
$L=\mathbb{Q}(\zeta_m)$. Then

$\mathrm{(i)}$ $\left(\frac{\alpha}{a}\right)_m=\left(\frac{a}{\alpha}\right)_m$,\qquad \qquad\;\; if\;\; $l>2$,

$\mathrm{(ii)}$ $\left(\frac{\alpha}{a}\right)_m=\left(\frac{(-1)^{(a-1)/2}a}{\alpha}\right)_m$, \;\;if\; \;$l=2$.
\end{thm}

\smallskip

\brem\label{re1}
$\mathrm{(i)}$ When $m=8$, let $\alpha\in\mathbb{Z}[\zeta_8]$  be coprime with $8$, we have
 $\alpha$ is primary if and only if $\alpha\equiv1$ or $1+\zeta_8+\zeta_8^3\pmod{2}$(see \cite[Th.14.2.1, p. 471]{ev}).

$\mathrm{(ii)}$ Let $m=2^n$ with $n\geq3$ and let $\alpha\in \mathbb{Z}[\zeta_m]$ be coprime with
$2$. There are exactly two $m$-th roots of unity $\mu$ for which $\mu\alpha$ is primary(see \cite[Th.14.6.2, p. 484]{ev}).

$\mathrm{(iii)}$ It is sufficient to apply Theorem \ref{Eis} to cases $m=8$ and $m=16$ in this paper.
\erem

\smallskip

\section{Explicit primality test}
From now on we will deduce an explicit primality test for $M=h\cdot 2^n\pm 1$ with $n\geq2$ and $h\nequiv 0$
$\pmod{17}$.
For any odd integer $k$ we set $k^*=(-1)^{(k-1)/2}k$. This notation allows us to treat
$h\cdot 2^n\pm 1$ simultaneously. And if $M=h\cdot 2^n\pm 1$, then $M^*=(\pm h)2^n+ 1$.

In this section let $\zeta_8 = e^{2\pi \sqrt{-1}/8}$ and $\zeta_{16}= e^{2\pi \sqrt{-1}/16}$,
and let $L_1 = \mathbb{Q}(\zeta_8)$ and $L_2=\mathbb{Q}(\zeta_{16})$
be the eighth and sixteenth cyclotomic fields respectively.
Let $D_1=\mathbb{Z}[\zeta_8]$ and $D_2=\mathbb{Z}[\zeta_{16}]$ be the corresponding cyclotomic rings.
$D_1$ and $D_2$ are both Principal Ideal Domains (PID) (see \cite[Th.11.1]{wa}).
Let $G=\mbox{Gal}(\mathbb{Q}(\zeta_{16})/\mathbb{Q})$ be the Galois
group of $\mathbb{Q}(\zeta_{16})$ over $\mathbb{Q}$. For every odd integer $c$ denote by $\sigma_c$ the
element of $G$ that sends $\zeta_{16}$ to $\zeta_{16}^c$. We also denote by $\sigma_c$ the element of Gal$(\mathbb{Q}(\zeta_{8})/\mathbb{Q})$ that sends $\zeta_{8}$ to $\zeta_{8}^c$. We know that Gal$(\mathbb{Q}(\zeta_{16})/\mathbb{Q})=\{\sigma_{\pm i}\mid i=1,3,5,7\}$ and Gal$(\mathbb{Q}(\zeta_{8})/\mathbb{Q})=\{\sigma_{\pm i}\mid i=1,3\}$. For $\tau$ in $\mathbb{Z}[G]$ and $\alpha$ in
$L_2$ with $\alpha\neq0$ we often denote by $\alpha^\tau$ to the action of the element $\tau$ of  $\mathbb{Z}[G]$ on the
element $\alpha$ of $L_2$, that is,
$$
\alpha^\tau:=\prod_{\sigma\in G}\sigma(\alpha)^{k_{\sigma}}, \text{ if }\tau=\sum_{\sigma\in G}k_{\sigma}\sigma\text{ where }k_{\sigma}\in\mathbb{Z}.
$$
If $\tau\in G$, we will either write $\alpha^\tau$ or $\tau(\alpha)$.
Since $L_1$ is contained in $L_2$, the element of $G$ can also act on elements of $L_1$. We also write $\sigma_1=1$ in $\mathbb{Z}[G]$.

Let $K_1=\mathbb{Q}(\zeta_8+\zeta_8^{-1})=\mathbb{Q}(\sqrt{2})$ and $K_2=\mathbb{Q}(\zeta_{16}+\zeta_{16}^{-1})$
be the maximal real subfield of $L_1$ and $L_2$ respectively. We know that Gal$(K_1/\mathbb{Q})=\{\sigma_i|_{K_1}\mid i=1,3\}$ and Gal$(K_2/\mathbb{Q})=\{\sigma_i|_{K_2}\mid i=1,3,5,7\}$.
Let $\pi_1\in D_1$ and $\pi_2\in D_2$ be two elements such that $\pi_1$, $\pi_2\notin \mathbb{R}$.
We denote two elements
$\alpha_1=(\pi_1/\bar{\pi_1})^{1+3\sigma_3}$ and
$\alpha_2=(\pi_2/\bar{\pi_2})^{1+3\sigma_{-5}+5\sigma_{-3}+7\sigma_7}$, where a bar indicates the complex conjugation.
Next we define some sequences.

$\mathrm{(1)}$ Sequences $\{T_k\}$ and $\{N_k\}$:
Set $T_0=Trace_{K_1/\mathbb{Q}}(\alpha_1^h+\bar{\alpha_1}^h)$ and $N_0=Norm_{K_1/\mathbb{Q}}(\alpha_1^h+\bar{\alpha_1}^h)$.
 For $k\geq0$ define $T_{k+1}$ and $N_{k+1}$
recursively by the formulas:
\beq \label{t} T_{k+1}=T_k^2-2N_k-4, \eeq
\beq \label{n} N_{k+1}=N_k^2-2T_k^2+4N_k+4.\eeq

$\mathrm{(2)}$ Sequences $\{X_k\}$, $\{Y_k\}$, $\{Z_k\}$ and $\{W_k\}$:
Set $X_0=Trace_{K_2/\mathbb{Q}}(\eta)$, $Y_0=\eta\sigma_3(\eta)+
\eta\sigma_5(\eta)+
\eta\sigma_7(\eta)
+\sigma_3(\eta)\sigma_5(\eta)+\sigma_3(\eta)\sigma_7(\eta)+\sigma_5(\eta)\sigma_7(\eta)$,
 $Z_0=\eta\sigma_3(\eta)\sigma_5(\eta)+\eta\sigma_3(\eta)\sigma_7(\eta)+\eta\sigma_5(\eta)\sigma_7(\eta)
 +\sigma_3(\eta)\sigma_5(\eta)\sigma_7(\eta)$ and $W_0=Norm_{K_2/\mathbb{Q}}(\eta)$, where $\eta=\alpha_2^h+\bar{\alpha_2}^h$. For $k\geq0$ define $X_{k+1}$, $Y_{k+1}$,
 $Z_{k+1}$ and $W_{k+1}$
recursively by the formulas:
\beq\label{x} X_{k+1}=X_k^2-2Y_k-8, \eeq
\beq\label{y} Y_{k+1}=Y_k^2-2X_kZ_k+2W_k-6X_k^2+12Y_k+24,\eeq
\beq\label{z} Z_{k+1}=Z_k^2-2W_kY_k-4Y_k^2+8X_kZ_k-8W_k+12X_k^2-24Y_k-32,\eeq
\beq\label{w} W_{k+1}=W_k^2-2Z_k^2+4W_kY_k+4Y_K^2-8X_kZ_k+8W_k-8X_k^2+16Y_k+16.\eeq

Our explicit primality test is described as follows:
\smallskip


\begin{thm}\label{Main}
Let $M=h\cdot 2^n\pm 1$ with $n\geq 7$, $0<h<2^{n-6}$, $h$ odd and $h \nequiv 0\pmod{17}$.
 Let $\pi_1 = 1+2\zeta_8^3$ and $\pi_2 =1-\zeta_{16}+\zeta_{16}^5$
in the above sequences $\{T_k\}$ and $\{N_k\}$, and sequences $\{X_k\}$,
 $\{Y_k\}$, $\{Z_k\}$ and $\{W_k\}$. Let $Q_i(i=1,\ldots,7)$ be seven integers satisfying $Q_i\equiv\pm5^{j2^{n-7}}\pmod{2^{n-3}}$ with $j=0,1,2,3$ and $1<Q_i<2^{n-3}$. Suppose that $M$ is not divisible by all $Q_i$ for $1\leq i \leq7$. Then $M$ is prime if and only if one of the following holds:

$\mathrm{(i)}$\;\; $M^* \equiv \pm 4\pmod{17}$, and $T_{n-3}\equiv -N_{n-3}\equiv -4\pmod{M}$.

$\mathrm{(ii)}$
\;$M^* \equiv \pm 2, \pm 8\pmod{17}$, and $T_{n-3}\equiv N_{n-3}\equiv 0\pmod{M}$ .

$\mathrm{(iii)}$
 $M^* \equiv \pm 3, \pm 5, \pm 6, \pm 7\pmod{17}$, and $T_{n-3}\equiv 0\pmod{M}$ and $N_{n-3}\equiv -2\pmod{M}$.

$\mathrm{(iv)}$
$M^* \equiv -1\pmod{17}$, and $X_{n-4}\equiv -8\pmod{M}$, $Y_{n-4}\equiv 24\pmod{M}$,
$Z_{n-4}\equiv -32\pmod{M}$ and $W_{n-4}\equiv 16\pmod{M}$.
\end{thm}

\smallskip

Before proving the theorem we first show some preliminary lemmas.

\smallskip


\blem\label{lemma1}
 Let $M=h\cdot 2^n\pm 1$ with $n\geq 3$, and let $\pi\in D_1$
with $gcd(\pi,2M)=1$. Set  $\alpha=(\pi/\bar{\pi})^{1+3\sigma_3}$.
Let $\{T_k\}$ and $\{N_k\}$  be the sequences defined in (\ref{t}) and (\ref{n}) with
$T_0=Trace_{K_1/\mathbb{Q}}(\alpha^h+\bar{\alpha}^h)$ and $N_0=Norm_{K_1/\mathbb{Q}}(\alpha^h+\bar{\alpha}^h)$.
Suppose $M$ is prime, then we have

$\mathrm{(i)}$ If $\left(\frac{\pi}{M}\right)_8=-1$, then $T_{n-3}\equiv -N_{n-3}\equiv -4\pmod{M}$.

$\mathrm{(ii)}$ If $\left(\frac{\pi}{M}\right)_8=\pm \zeta_8^2$, then $T_{n-3}\equiv N_{n-3}\equiv 0\pmod{M}$.

$\mathrm{(iii)}$ If $\left(\frac{\pi}{M}\right)_8=\pm\zeta_8,\pm\zeta_8^{-1}$, then $T_{n-3}\equiv 0\pmod{M}$ and $N_{n-3}\equiv -2\pmod{M}$.
\elem


 \bpf{}\; Since $M$ is a prime, when $M=h\cdot 2^n + 1\equiv 1\pmod{8}$, the ideal $MD_1$ factors in $D_1$ as a product of 4 distinct prime ideals. We write
$MD_1=(\mathfrak{p}\bar{\mathfrak{p}})^{1+\sigma_3}$, thus

\[
\begin{array}{ll}
\left(\frac{\pi}{M}\right)_8
&=\left(\frac{\pi}{(\mathfrak{p}\bar{\mathfrak{p}})^{1+\sigma_3}}\right)_8
=\left(\frac{\pi}{\mathfrak{p}\bar{\mathfrak{p}}(\mathfrak{p}\bar{\mathfrak{p}})^{\sigma_3}}\right)_8\\\\
&=\left(\frac{\pi/\bar{\pi}}{\mathfrak{p}}\right)_8 \left(\frac{(\pi/\bar{\pi})^{3\sigma_3}}{\mathfrak{p}}\right)_8\\\\
&=\left(\frac{\alpha}{\mathfrak{p}}\right)_8\equiv\alpha^{\frac{M-1}{8}}
\equiv\alpha^{h\cdot2^{n-3}} \pmod{\mathfrak{p}}.
\end{array}
\]

Since $\mathfrak{p}$ is an arbitrary prime ideal lying over $M$, we have
$$\left(\frac{\pi}{M}\right)_8\equiv\alpha^{h\cdot2^{n-3}} \pmod{M}.$$

When $M=h\cdot 2^n - 1\equiv -1\pmod{8}$, the ideal $MD_1$ factors in $D_1$ as
a product of 2 distinct prime ideals. Write
$MD_1=\mathfrak{p}\mathfrak{p}^{\sigma_3}$, then

\[
\begin{array}{ll}
\left(\frac{\pi}{M}\right)_8
&=\left(\frac{\pi}{\mathfrak{p}\mathfrak{p}^{\sigma_3}}\right)_8=
\left(\frac{\pi}{\mathfrak{p}}\right)_8 \left(\frac{\pi}{\mathfrak{p}^{\sigma_3}}\right)_8\\\\
&=\left(\frac{\pi^{1+3\sigma_3}}{\mathfrak{p}}\right)_8\equiv(\pi^{1+3\sigma_3})^{\frac{M^2-1}{8}}\\\\
&\equiv\bar{\alpha}^{\frac{M+1}{8}}\equiv\bar{\alpha}^{h\cdot2^{n-3}} \pmod{\mathfrak{p}}.
\end{array}
\]

The last second congruence holds because of $\pi^M\equiv\bar{\pi}\pmod{\mathfrak{p}}$,
it can be seen by observing that the complex conjugation coincides with the
 Frobenius automorphism of $D_1/\mathfrak{p}$. As before  we obtain
$$\left(\frac{\pi}{M}\right)_8\equiv\bar{\alpha}^{h\cdot2^{n-3}} \pmod{M}.$$
Hence for $M=h\cdot 2^n\pm 1$, we always have
\beq\label{8}\alpha^{h2^{n-3}}+\bar{\alpha}^{h2^{n-3}}\equiv\left(\frac{\pi}{M}\right)_8+\left(\frac{\pi}{M}\right)_8^{-1} \pmod{M}.\eeq

For $k\geq 0$ let $T_k=Trace_{K_1/\mathbb{Q}}(\alpha^{h2^k}+\bar{\alpha}^{h2^k})$ and
$N_k=Norm_{K_1/\mathbb{Q}}(\alpha^{h2^k}+\bar{\alpha}^{h2^k})$. We claim that $T_k$
and $N_k$ satisfy the recurrent relations given by (\ref{t}) and (\ref{n}). To see this
we let $A_k=\alpha^{h2^k}+\bar{\alpha}^{h2^k}$ and $B_k=\sigma_3(A_k)$. Thus
$T_k=A_k+B_k$ and $N_k=A_kB_k$.

By computation, we have $A_{k+1}=A_k^2-2$ and $B_{k+1}=B_k^2-2$. We substitute these in $T_{k+1}$ and
$N_{k+1}$, and obtain
$$T_{k+1}=A_k^2+B_k^2-4=T_k^2-2N_k-4,$$
$$N_{k+1}=(A_k^2-2)(B_k^2-2)=N_k^2-2(T_k^2-2N_k)+4=N_k^2-2T_k^2+4N_k+4.$$
Since we have proved that $T_k$ and $N_k$ satisfy the recurrence relations given by
(\ref{t}) and (\ref{n}), (\ref{8}) implies that
$$T_{n-3}\equiv\left[\left(\frac{\pi}{M}\right)_8+\left(\frac{\pi}{M}\right)_8^{-1}\right]+\left[\left(\frac{\pi}{M}\right)_8^3+
\left(\frac{\pi}{M}\right)_8^{-3}\right]\pmod{M},$$
$$N_{n-3}\equiv\left[\left(\frac{\pi}{M}\right)_8+\left(\frac{\pi}{M}\right)_8^{-1}\right]\cdot\left[\left(\frac{\pi}{M}\right)_8^3+
\left(\frac{\pi}{M}\right)_8^{-3}\right]\pmod{M}.$$
Hence we get if $\left(\frac{\pi}{M}\right)_8=-1$, then $T_{n-3}\equiv-4\pmod{M}$ and $N_{n-3}\equiv4\pmod{M}$.
If $\left(\frac{\pi}{M}\right)_8=\pm \zeta_8^2$, then $T_{n-3}\equiv0\pmod{M}$ and $N_{n-3}\equiv0\pmod{M}$.
If $\left(\frac{\pi}{M}\right)_8=\pm\zeta_8,\pm\zeta_8^{-1}$, then
$T_{n-3}\equiv0\pmod{M}$ and $N_{n-3}\equiv-2\pmod{M}$. This completes the proof of
the three cases.\qed

\epf
\smallskip


\blem\label{lemma2}
Let $M=h\cdot 2^n\pm 1$ with $n\geq 4$, and let $\pi\in D_2$ with $gcd(\pi,2M)=1$.
Set $\alpha=(\pi/\bar{\pi})^{1+3\sigma_{-5}+5\sigma_{-3}+7\sigma_7}$.
Let $\{X_k\}$, $\{Y_k\}$, $\{Z_k\}$ and $\{W_k\}$  be the sequences defined in (\ref{x}), (\ref{y}), (\ref{z}) and (\ref{w}) with $X_0=Trace_{K_2/\mathbb{Q}}(\eta)$, $Y_0=\eta\sigma_3(\eta)+
\eta\sigma_5(\eta)+\eta\sigma_7(\eta)
+\sigma_3(\eta)\sigma_5(\eta)+\sigma_3(\eta)\sigma_7(\eta)+\sigma_5(\eta)\sigma_7(\eta)$,
$Z_0=\eta\sigma_3(\eta)\sigma_5(\eta)+\eta\sigma_3(\eta)\sigma_7(\eta)+\eta\sigma_5(\eta)\sigma_7(\eta)
+\sigma_3(\eta)\sigma_5(\eta)\sigma_7(\eta)$ and $W_0=Norm_{K_2/\mathbb{Q}}(\eta)$,
where $\eta=\alpha^h+\bar{\alpha}^h$.
Suppose $M$ is prime and $\left(\frac{\pi}{M}\right)_{16}=-1$, then  $X_{n-4}\equiv -8\pmod{M}$, $Y_{n-4}\equiv 24\pmod{M}$,
$Z_{n-4}\equiv -32\pmod{M}$ and $W_{n-4}\equiv 16\pmod{M}$.

 \elem


 \bpf{}\;  Since $M$ is a prime, when $M=h\cdot 2^n + 1\equiv 1\pmod{16}$, the ideal $MD_2$ factors in $D_2$ as a product of 8 distinct prime ideals. We write
$MD_2=(\mathfrak{p}\bar{\mathfrak{p}})^{1+\sigma_3+\sigma_5+\sigma_7}$, thus

\[
\begin{array}{ll}
\left(\frac{\pi}{M}\right)_{16}
&=\left(\frac{\pi}{(\mathfrak{p}\bar{\mathfrak{p}})^{1+\sigma_3+\sigma_5+\sigma_7}}\right)_{16}\\\\
&=\left(\frac{(\pi/\bar{\pi})^{1+3\sigma_{-5}+5\sigma_{-3}+7\sigma_7}}{\mathfrak{p}}\right)_{16}\\\\
&=\left(\frac{\alpha}{\mathfrak{p}}\right)_{16}\equiv\alpha^{\frac{M-1}{16}}
\equiv\alpha^{h\cdot2^{n-4}} \pmod{\mathfrak{p}}.
\end{array}
\]

Since $\mathfrak{p}$ is an arbitrary prime ideal lying over $M$, we have
$$\left(\frac{\pi}{M}\right)_{16}\equiv\alpha^{h\cdot2^{n-4}} \pmod{M}.$$

When $M=h\cdot 2^n - 1\equiv -1\pmod{16}$, the ideal $MD_2$ factors in $D_2$ as
a product of 4 distinct prime ideals. We write
$MD_2=\mathfrak{p}^{1+\sigma_3+\sigma_5+\sigma_7}$, hence

\[
\begin{array}{ll}
\left(\frac{\pi}{M}\right)_{16}
&=\left(\frac{\pi}{\mathfrak{p}}\right)_{16} \left(\frac{\pi}{\mathfrak{p}^{\sigma_3}}\right)_{16}
\left(\frac{\pi}{\mathfrak{p}^{\sigma_5}}\right)_{16}\left(\frac{\pi}{\mathfrak{p}^{\sigma_7}}\right)_{16}\\\\
&=\left(\frac{\pi^{1+3\sigma_{-5}+5\sigma_{-3}+7\sigma_7}}{\mathfrak{p}}\right)_{16}
\equiv\left(\pi^{1+3\sigma_{-5}+5\sigma_{-3}+7\sigma_7}\right)^{(M^2-1)/16}\\\\
&\equiv\bar{\alpha}^{(M+1)/16}\equiv\bar{\alpha}^{h\cdot2^{n-4}} \pmod{\mathfrak{p}}.
\end{array}
\]

The last second congruence holds because of $\pi^M\equiv\bar{\pi}\pmod{\mathfrak{p}}$,
it can be seen by the same observation as in Lemma \ref{lemma1}. And  we obtain
$$\left(\frac{\pi}{M}\right)_{16}\equiv\bar{\alpha}^{h\cdot2^{n-4}} \pmod{M}.$$
Hence for $M=h\cdot 2^n\pm 1$, we always have
\beq\label{b}\alpha^{h2^{n-4}}+\bar{\alpha}^{h2^{n-4}}\equiv\left(\frac{\pi}{M}\right)_{16}+\left(\frac{\pi}{M}\right)_{16}^{-1} \pmod{M}.\eeq

For $k\geq 0$ let $X_k=Trace_{K_2/\mathbb{Q}}(\eta_k)$,
$Y_k=\eta_k\sigma_3(\eta_k)+
\eta_k\sigma_5(\eta_k)+\eta_k\sigma_7(\eta_k)
+\sigma_3(\eta_k)\sigma_5(\eta_k)+\sigma_3(\eta_k)\sigma_7(\eta_k)+\sigma_5(\eta_k)\sigma_7(\eta_k)$,
$Z_k=\eta_k\sigma_3(\eta_k)\sigma_5(\eta_k)+\eta_k\sigma_3(\eta_k)\sigma_7(\eta_k)+\eta_k\sigma_5(\eta_k)\sigma_7(\eta_k)
+\sigma_3(\eta_k)\sigma_5(\eta_k)\sigma_7(\eta_k)$ and
$W_k=Norm_{K_2/\mathbb{Q}}(\eta_k)$, where
$\eta_k=\alpha^{h2^k}+\bar{\alpha}^{h2^k}\in K_2$.

We claim that $X_k$, $Y_k$, $Z_k$
and $W_k$ satisfy the recurrent relations given by (\ref{x}), (\ref{y}), (\ref{z}) and (\ref{w}). To see this
we let $A_k=\eta_k$, $B_k=\sigma_3(\eta_k)$, $C_k=\sigma_5(\eta_k)$ and $D_k=\sigma_7(\eta_k)$. So
$X_k=A_k+B_k+C_k+D_k$, $Y_k=A_kB_k+A_kC_k+A_kD_k+B_kC_k+B_kD_k+C_kD_k$,
 $Z_k=A_kB_kC_k+A_kB_kD_k+A_kC_kD_k+B_kC_kD_k$ and $W_k=A_kB_kC_kD_k$.

By computation, we get $A_{k+1}=A_k^2-2$, $B_{k+1}=B_k^2-2$, $C_{k+1}=C_k^2-2$ and $D_{k+1}=D_k^2-2$,
substituting these in $X_{k+1}$, $Y_{k+1}$, $Z_{k+1}$ and
$W_{k+1}$, then
\[
\begin{array}{ll}
X_{k+1}
&=A_k^2+B_k^2+C_k^2+D_k^2-8\\\\
&=X_k^2-2Y_k-8,
\end{array}
\]
\[
\begin{array}{ll}
Y_{k+1}
&=(A_kB_k)^2+(A_kC_k)^2+(A_kD_k)^2+(B_kC_k)^2+(B_kD_k)^2+(C_kD_k)^2\\\\&-6(A_k^2+B_k^2+C_k^2+D_k^2)+24\\\\
&=Y_k^2-2(X_kZ_k-W_k)-6(X_k^2-2Y_k)+24,
\end{array}
\]

\[
\begin{array}{ll}
Z_{k+1}
&=(A_kB_kC_k)^2+(A_kB_kD_k)^2+(A_kC_kD_k)^2+(B_kC_kD_k)^2\\\\&-4[(A_kB_k)^2+(A_kC_k)^2+(A_kD_k)^2
+(B_kC_k)^2+(B_kD_k)^2+(C_kD_k)^2]\\\\&+12(A_k^2+B_k^2+C_k^2+D_k^2)-32\\\\
&=Z_k^2-2W_kY_k-4Y_k^2+8(X_kZ_k-W_k)+12X_k^2-24Y_k-32,
\end{array}
\]

\[
\begin{array}{ll}
W_{k+1}
&=(A_k^2-2)(B_k^2-2)(C_k^2-2)(D_k^2-2)\\\\
&=W_k^2-2Z_k^2+4W_kY_k+4Y_k^2-8(X_kZ_k-W_k)-8X_k^2+16Y_k+16.
\end{array}
\]

After taking $\left(\frac{\pi}{M}\right)_{16}=-1$ in (\ref{b}) we have
\beq\label{16}\alpha^{h2^{n-4}}+\bar{\alpha}^{h2^{n-4}}\equiv -2\pmod{M}.\eeq
Since we have proved that $X_k$, $Y_k$, $Z_k$
and $W_k$ satisfy the recurrent relations given by (\ref{x}), (\ref{y}), (\ref{z}) and (\ref{w}).
Now substituting  (\ref{16}) in $X_{n-4}$, $Y_{n-4}$, $Z_{n-4}$ and $W_{n-4}$,
 we get
$$X_{n-4}\equiv-2\cdot4\equiv-8\pmod{M},\quad Y_{n-4}\equiv4\cdot6\equiv 24 \pmod{M},$$
$$Z_{n-4}\equiv-8\cdot4\equiv -32\pmod{M}, \quad
W_{n-4}\equiv(-2)^4\equiv 16\pmod{M}.$$\qed

\epf
\smallskip


\blem\label{lemma3}
Let the number fields $L_1$ and $K_1$ be as before, let $q$ be an odd rational prime
and let $\pi\in D_1$ be prime to $q$. Set $\alpha=\pi/\bar{\pi}$. Let $\{T_k\}$ and
$\{N_k\}$ be the sequences defined in (\ref{t}) and (\ref{n}) with
$T_0=Trace_{K_1/\mathbb{Q}}(\alpha+\bar{\alpha})$ and $N_0=Norm_{K_1/\mathbb{Q}}(\alpha+\bar{\alpha})$.
Suppose that for some $j\geq0$, one of the following statements holds:

$\mathrm{(i)}$  $T_j\equiv -N_j\equiv -4\pmod{q}$,

$\mathrm{(ii)}$  $T_j\equiv N_j\equiv 0\pmod{q}$,

$\mathrm{(iii)}$ $T_j\equiv 0\pmod{q}$ and $N_j\equiv -2\pmod{q}$.

Then $q^2\equiv1\pmod{2^{j+1}}$.

\elem


\bpf{}\;
By Lemma \ref{lemma1}, we have $T_j=Trace_{K_1/\mathbb{Q}}(\alpha^{2^j}+\bar{\alpha}^{2^j})$ and
$N_j=Norm_{K_1/\mathbb{Q}}(\alpha^{2^j}+\bar{\alpha}^{2^j})$.
 Let $\mathfrak{q}$ be a prime ideal in
the ring of integers of $K_1$ lying over $q$, and $\mathfrak{Q}$ be a prime ideal
of $D_1$ lying over $\mathfrak{q}$. Let $\beta=\alpha^{2^j}+\bar{\alpha}^{2^j}$,
then we have

$\mathrm{(i)}$ $T_j\equiv -N_j\equiv -4\pmod{q}$ means
$Trace_{K_1/\mathbb{Q}}(\beta)\equiv -4\equiv -Norm_{K_1/\mathbb{Q}}(\beta)\pmod{\mathfrak{q}}$, which implies
$\beta^2+4\beta+4\equiv0\pmod{\mathfrak{q}}$, i.e., $\alpha^{2^j}+\bar{\alpha}^{2^j}\equiv-2\pmod{\mathfrak{Q}}$. Multiplying both sides
of the congruence by $\alpha^{2^j}=\bar{\alpha}^{-2^j}$ gives
$$\alpha^{2^j}\equiv-1\pmod{\mathfrak{Q}}.$$
It implies that the image of $\alpha$ has order $2^{j+1}$ in the multiplicative group $(D_1/\mathfrak{Q})^*$.
This group has order N$(\mathfrak{Q})-1$ which divides
$q^2-1$, i.e., $q^2\equiv1\pmod{2^{j+1}}$.

$\mathrm{(ii)}$ $T_j\equiv N_j\equiv 0\pmod{q}$ means
$Trace_{K_1/\mathbb{Q}}(\beta)\equiv 0\equiv Norm_{K_1/\mathbb{Q}}(\beta)\pmod{\mathfrak{q}}$, which implies
$\beta^2\equiv0\pmod{\mathfrak{q}}$, i.e., $\alpha^{2^j}+\bar{\alpha}^{2^j}\equiv0\pmod{\mathfrak{Q}}$. Also multiplying both sides
of the congruence by $\alpha^{2^j}=\bar{\alpha}^{-2^j}$ gives
$$\alpha^{2^{j+1}}\equiv-1\pmod{\mathfrak{Q}}.$$
It implies that the image of $\alpha$ has order $2^{j+2}$ in  $(D_1/\mathfrak{Q})^*$. Hence $2^{j+2}$ must divides
$q^2-1$, i.e., $q^2\equiv1\pmod{2^{j+2}}$.

$\mathrm{(iii)}$ $T_j\equiv 0\pmod{q}$ and $N_j\equiv -2\pmod{q}$ imply
$Trace_{K_1/\mathbb{Q}}(\beta)\equiv 0\pmod{\mathfrak{q}}$ and $Norm_{K_1/\mathbb{Q}}(\beta)\equiv -2\pmod{\mathfrak{q}}$, which deduce
$\beta^2-2\equiv0\pmod{\mathfrak{q}}$, i.e., $\alpha^{2^{j+1}}+\bar{\alpha}^{2^{j+1}}\equiv0\pmod{\mathfrak{Q}}$. Also we obtain
$$\alpha^{2^{j+2}}\equiv-1\pmod{\mathfrak{Q}}.$$
That is the image of $\alpha$ has order $2^{j+3}$ in group $(D_1/\mathfrak{Q})^*$.
And we reach $q^2\equiv1\pmod{2^{j+3}}$. This completes the proof.\qed

\epf
\smallskip


\blem\label{lemma4}
Let the number fields $L_2$ and $K_2$ be as before, let $q$ be an odd rational prime
and let $\pi\in D_2$ be prime to $q$. Set $\alpha=\pi/\bar{\pi}$. Let $\{X_k\}$,
$\{Y_k\}$, $\{Z_k\}$ and $\{W_k\}$  be the sequences defined in (\ref{x}), (\ref{y}), (\ref{z}) and (\ref{w}) with
 $X_0=Trace_{K_2/\mathbb{Q}}(\eta)$, $Y_0=\eta\sigma_3(\eta)+
\eta\sigma_5(\eta)+\eta\sigma_7(\eta)
+\sigma_3(\eta)\sigma_5(\eta)+\sigma_3(\eta)\sigma_7(\eta)+\sigma_5(\eta)\sigma_7(\eta)$,
$Z_0=\eta\sigma_3(\eta)\sigma_5(\eta)+\eta\sigma_3(\eta)\sigma_7(\eta)+\eta\sigma_5(\eta)\sigma_7(\eta)
+\sigma_3(\eta)\sigma_5(\eta)\sigma_7(\eta)$ and $W_0=Norm_{K_2/\mathbb{Q}}(\eta)$,
where $\eta=\alpha+\bar{\alpha}$.
Suppose that for some $j\geq0$, $X_j\equiv -8\pmod{q}$, $Y_j\equiv 24\pmod{q}$,
$Z_j\equiv -32\pmod{q}$, and $W_j\equiv 16\pmod{q}$. Then $q^4\equiv1\pmod{2^{j+1}}$.

\elem


\bpf{}\;
By Lemma \ref{lemma2}, we have $X_j=Trace_{K_2/\mathbb{Q}}(\beta)$,
$Y_j=\beta\sigma_3(\beta)+
\beta\sigma_5(\beta)+\beta\sigma_7(\beta)
+\sigma_3(\beta)\sigma_5(\beta)+\sigma_3(\beta)\sigma_7(\beta)+\sigma_5(\beta)\sigma_7(\beta)$,
$Z_j=\beta\sigma_3(\beta)\sigma_5(\beta)+\beta\sigma_3(\beta)\sigma_7(\beta)+
\beta\sigma_5(\beta)\sigma_7(\beta)
+\sigma_3(\beta)\sigma_5(\beta)\sigma_7(\beta)$ and
$W_j=Norm_{K_2/\mathbb{Q}}(\beta)$, where
$\beta=\alpha^{2^j}+\bar{\alpha}^{2^j}\in K_2$.
Let $\mathfrak{q}$ be a prime ideal in
the ring of integers of $K_2$ lying over $q$, and $\mathfrak{Q}$ be a prime ideal
of $D_2$ lying over $\mathfrak{q}$.

By the assumption we get
$Trace_{K_2/\mathbb{Q}}(\beta)\equiv -8\pmod{q}$,
$\beta\sigma_3(\beta)+
\beta\sigma_5(\beta)+\beta\sigma_7(\beta)
+\sigma_3(\beta)\sigma_5(\beta)+\sigma_3(\beta)\sigma_7(\beta)+
\sigma_5(\beta)\sigma_7(\beta)\equiv 24\pmod{q}$,
$\beta\sigma_3(\beta)\sigma_5(\beta)+\beta\sigma_3(\beta)\sigma_7(\beta)+
\beta\sigma_5(\beta)\sigma_7(\beta)
+\sigma_3(\beta)\sigma_5(\beta)
\sigma_7(\beta)\equiv -32\pmod{q}$ and
$Norm_{K_2/\mathbb{Q}}(\beta)\equiv 16\pmod{q}$, which implies
$(\beta+2)^4=\beta^4+8\beta^3+24\beta^2+32\beta+16 \equiv0\pmod
{\mathfrak{q}}$, i.e., $\alpha^{2^j}+\bar{\alpha}^{2^j}\equiv-2\pmod{\mathfrak{Q}}$. Multiplying both sides
of the congruence by $\alpha^{2^j}=\bar{\alpha}^{-2^j}$ gives
$$\alpha^{2^j}\equiv-1\pmod{\mathfrak{Q}}.$$
It implies that the image of $\alpha$ has order $2^{j+1}$ in the multiplicative group $(D_2/\mathfrak{Q})^*$.
The order of this group is N$(\mathfrak{Q})-1$ which divides
$q^4-1$, i.e., $q^4\equiv1\pmod{2^{j+1}}$.\qed

\epf
\smallskip

\brem
To prove the main theorem, we work in $L_1=\mathbb{Q}(\zeta_8)$ and $L_2=\mathbb{Q}(\zeta_{16})$.
When $\pi_1=1+2\zeta_8^3$ and $\pi_2=1-\zeta_{16}+\zeta_{16}^5$,
note that $Norm_{L_1/\mathbb{Q}}(\pi_1)=Norm_{L_2/\mathbb{Q}}(\pi_2)=17$.
By Remark \ref{re1}, $\pi_1\equiv1\pmod{2}$
implies that $\pi_1$ is a primary prime in $D_1$.
Since $2-\zeta_8=-\pi_1\cdot\zeta_8$, then
 $2\equiv\zeta_8\pmod{\pi_1}$. The verification of a primary element is quite troublesome.
Actually for $\pi_2$ we do not need to know whether it is primary or not. We can see this
from the process of the next proof. The choice of $\pi_2$ is enough for our explicit primality test.
\erem

\smallskip


\bpf{ \it (of Theorem \ref{Main})}\;
We first show that the congruences of the sequences are necessary for primality of
$M$. Suppose then that $M$ is a prime. Since $n\geq 7$, we have $M\neq17$, so
the hypotheses allow $M^*\equiv-1, \pm2, \pm3, \pm4, \pm5, \pm6, \pm7, \pm8\pmod
{17}$, hence mod $\pi_1$ and mod $\pi_2$.

By Remark \ref{re1}, there is a $16$-th root of unity $\mu$ such that $\mu\pi_2$ is
a primary element in $D_2$. Let $\pi=\mu\pi_2$. Since $\pi_1$ and $\pi$ are
primary primes in $D_1$ and $D_2$ respectively, we can apply Theorem \ref{Eis}(Eisenstein's
Reciprocity Law) to obtain
$$\left(\frac{\pi_1}{M}\right)_8=\left(\frac{M^*}{\pi_1}\right)_8 \quad
\mbox{and} \quad
\left(\frac{\pi}{M}\right)_{16}=\left(\frac{M^*}{\pi}\right)_{16}.$$
Note N$(\pi_1D_1)=\mbox{N}(\pi D_2)=17$,
now we compute as follows:

$\mathrm{(i)}$ Suppose $M^*\equiv\pm4\pmod{17}$, then
$\left(\frac{M^*}{\pi_1}\right)_8\equiv (M^*)^{(17-1)/8}\equiv(M^*)^2\equiv-1\pmod{\pi_1}$. And $\left(\frac{\pi_1}{M}\right)_8
=\left(\frac{M^*}{\pi_1}\right)_8=-1$, applying Lemma \ref{lemma1}$\mathrm{(i)}$ we get $T_{n-3}\equiv-4\equiv-N_{n-3}\pmod{M}$.

$\mathrm{(ii)}$ Suppose $M^*\equiv\pm2, \pm8\pmod{17}$,
 then $\left(\frac{M^*}{\pi_1}\right)_8\equiv(M^*)^2\equiv\pm4\equiv\pm\zeta_8^2\pmod{\pi_1}$. And $\left(\frac{\pi_1}{M}\right)_8
=\left(\frac{M^*}{\pi_1}\right)_8=\pm\zeta_8^2$, hence by Lemma \ref{lemma1}$\mathrm{(ii)}$ we get $T_{n-3}\equiv 0\equiv N_{n-3}\pmod{M}$.

$\mathrm{(iii)}$ Suppose $M^*\equiv\pm3, \pm5, \pm6, \pm7\pmod{17}$,
then $\left(\frac{M^*}{\pi_1}\right)_8\equiv(M^*)^2\equiv\pm2, \pm8\equiv\pm\zeta_8, \pm\zeta_8^{-1}\pmod{\pi_1}$. Thus $\left(\frac{\pi_1}{M}\right)_8
=\left(\frac{M^*}{\pi_1}\right)_8=\pm\zeta_8, \pm\zeta_8^{-1}$, also by Lemma \ref{lemma1}$\mathrm{(iii)}$ we obtain
$T_{n-3}\equiv 0\pmod{M}$ and $N_{n-3}\equiv -2\pmod{M}$.

$\mathrm{(iv)}$ Suppose $M^*\equiv-1\pmod{17}$, then
$\left(\frac{M^*}{\pi}\right)_{16}\equiv(M^*)^{(17-1)/16}\equiv M^*\equiv-1\pmod{\pi}$. So $\left(\frac{\pi}{M}\right)_{16}
=\left(\frac{M^*}{\pi}\right)_{16}=-1$. By the proof of Lemma \ref{lemma2}, we have
$$\alpha^{h2^{n-4}}+\bar{\alpha}^{h2^{n-4}}\equiv\left(\frac{\pi}{M}\right)_{16}+\left(\frac{\pi}{M}\right)_{16}^{-1}=-2 \pmod{M},$$
where $\alpha=\left(\frac{\pi}{\bar{\pi}}\right)^{1+3\sigma_{-5}+5\sigma_{-3}+7\sigma_7}$.

Let $\alpha_2=\left(\frac{\pi_2}{\bar{\pi_2}}\right)^{1+3\sigma_{-5}+5\sigma_{-3}+7\sigma_7}$, by computation,
\[
\begin{array}{ll}
\alpha^{h2^{n-4}}
&=\left(\frac{\pi}{\bar{\pi}}\right)^{2^{n-4}h(1+3\sigma_{-5}+5\sigma_{-3}+7\sigma_7)}\\\\
&=\left(\frac{\mu}{\bar{\mu}}\right)^{2^{n-4}h(1+3\sigma_{-5}+5\sigma_{-3}+7\sigma_7)}
\left(\frac{\pi_2}{\bar{\pi_2}}\right)^{2^{n-4}h(1+3\sigma_{-5}+5\sigma_{-3}+7\sigma_7)}\\\\
&=\mu^{2^{n-3}h(1+3\sigma_{-5}+5\sigma_{-3}+7\sigma_7)}\alpha_2^{h2^{n-4}}\\\\
&=\alpha_2^{h2^{n-4}}.
\end{array}
\]

The last second equality holds because $\mu$ is a $16$-th root of unity and $\bar{\mu}=\mu^{-1}$.
The last equality holds because of $n\geq7$ and $\mu^{2^{n-3}}=1$. Hence we have
$$\alpha_2^{h2^{n-4}}+\bar{\alpha_2}^{h2^{n-4}}\equiv-2 \pmod{M}.$$
And again by the proof of Lemma \ref{lemma2}, we obtain
$X_{n-4}\equiv -8\pmod{M}$, $Y_{n-4}\equiv 24\pmod{M}$,
$Z_{n-4}\equiv -32\pmod{M}$ and $W_{n-4}\equiv 16\pmod{M}$.
This completes the proof of necessity.

We now turn to the proof of sufficiency. Let $q$ be an arbitrary prime divisor of $M$.
In the first three cases, the hypotheses imply $q$ prime to $17$. Then take
 $\alpha=(\pi_1/\bar{\pi_1})^{h(1+3\sigma_3)}$ in Lemma \ref{lemma3}, we get  $q^2\equiv1\pmod{2^{n-2}}$. In the last case, let $\alpha=(\pi_2/\bar{\pi_2})^{h(1+3\sigma_{-5}+5\sigma_{-3}+7\sigma_7)}$ in Lemma
\ref{lemma4}, we obtain  $q^4\equiv1\pmod{2^{n-3}}$. Among all cases we always have
 $q^4\equiv1\pmod{2^{n-3}}$. By the assumption $M$ is not divisible by all
 $Q_i$ for $1\leq i\leq7$, which are all solutions of equation $x^4\equiv1\pmod{2^{n-3}}$
 between $1$ and $2^{n-3}$.  Then $q\geq2^{n-3}+1$ and $q^2\geq 2^{2n-6}+2^{n-2}+1
 =2^{n}(2^{n-6}+\frac{1}{4})+1>h\cdot 2^n+1\geq M$.
 Thus $q >\sqrt{M}$ for arbitrary prime
 divisor $q$ of $M$, that is to say $M$ is prime. This completes the proof of sufficiency.\qed

\epf

\vspace{5mm}

\noindent \textbf{Acknowledgments}\quad  The work of this paper
was supported by the NNSF of China (Grants Nos. 11071285,
61121062), 973 Project (2011CB302401) and the National Center for Mathematics and Interdisciplinary Sciences, CAS.


\end{document}